\documentclass{article}
\usepackage{amssymb}
\usepackage{amsmath}
\usepackage{amsfonts, latexsym}

\usepackage{fullpage}
\usepackage[utf8]{inputenc}

\usepackage{tikz-cd}
\usepackage{graphicx}
\usepackage{parskip}
\usepackage{enumitem}
\usepackage{verbatim,enumerate}
\usepackage{listings}
\usepackage{color}
\usepackage{hyperref}
\usepackage{cancel}

\usepackage{framed}
\usepackage{relsize}
\usepackage[]{algorithm2e}
\usepackage{amsmath,amssymb,latexsym,amsthm}
\hypersetup{final}
\usepackage[T1]{fontenc}
\usepackage{enumitem}
\usepackage{hyperref} 

\DeclareSymbolFont{symbolsC}{U}{txsyc}{m}{n}
\DeclareMathSymbol{\notniFromTxfonts}{\mathrel}{symbolsC}{61}

\newtheorem{thm}{Theorem}

\newtheorem{dfn}[thm]{Definition}
\newtheorem{exm}{Example}

\newcommand{\R}{\mathbb{R}}

\newcommand{\bea}{\begin{eqnarray}}
\newcommand{\eea}{\end{eqnarray}}
\newcommand{\bean}{\begin{eqnarray*}}
\newcommand{\eean}{\end{eqnarray*}}
\newcommand{\beq}{\begin{equation}}
\newcommand{\eeq}{\end{equation}}
\newcommand{\bac}{\begin{array}{c}}
\newcommand{\ball}{\begin{array}{ll}}
\newcommand{\ea}{\end{array}}

\renewcommand{\theequation}{\thesection.\arabic{equation}}

\renewcommand\qedsymbol{$\blacksquare$}

\newcommand{\bbR}{{\mathbb R}}

\newcommand{\eps}{\varepsilon}
\def\({\left(}
\def\){\right)}

\def\bl{\left\{}
\def\br{\right\}}
\def\ml{\left|}
\def\mr{\right|}

\def\ep{\varepsilon}

\usepackage{graphicx} 
\usepackage{mathrsfs}
\usepackage{bbm}
\usepackage{exscale}
\usepackage{hyperref}
\usepackage{animate}

\usepackage{xcolor}

\theoremstyle{Definition}
\newtheorem{example}{Example}
\newtheorem{problem}{Problem}

\def\ep{\varepsilon}

\def\({\left(}
\def\){\right)}

\def\bl{\left\{}
\def\br{\right\}}
\def\ml{\left|}
\def\mr{\right|}

\def\ep{\varepsilon}

\def\eps{\varepsilon}

\begin{document}

\title{On Whitney extensions and orthonormal expansions}
\date{\today}
\author{Steven B. Damelin \thanks{Department of Mathematics, University of Michigan, 530 Church Street, Ann Arbor, MI, USA.\, email: damelin@umich.edu.}\thanks{
}}

\maketitle

\thispagestyle{empty}
\parskip=10pt

\begin{abstract}

The Whitney near extension problem for finite sets in $\mathbb R^d,\, d\geq 2$ asks the following: Let $\phi:E\to \mathbb R^d$ be a near distortion on a finite set $E\subset \mathbb R^d$ with certain geometry.
How to decide whether $\phi$ extends to a smooth, one to one and onto near distortion $\Phi:\mathbb R^d\to \mathbb R^d$ which agrees with $\phi$ on $E$ and with Euclidean motions in $\mathbb R^d$. The Whitney near extension problem for compact sets $E\subset U$ in open subsets $U$ of $\mathbb R^n,\, n\geq 1$ asks the following: 
Let $U\subset R^n$ be open and let $E\subset U$ be a compact set. Let $\phi:U\to \mathbb R^n$ be a smooth near isometry. How to decide if there exists a smooth one-to-one and onto 
near isometry $\Phi:\mathbb R^n\to \mathbb R^n$ which extends $\phi$ on $E$ and agrees with Euclidean motions on $\mathbb R^n$. The classical Whitney extension problem asks the following: Let $\phi:E\to \mathbb R$ be a map defined on an arbitrary set $E\subset \mathbb R^n$. How can one decide whether $\phi$ extends to a map $\Phi:\mathbb R^n\to \mathbb R$ which agrees with $\phi$ on $E$ and is in $C^m(\mathbb R^n),\, m\geq 1$, the space of functions from $\mathbb R^n$ to $\mathbb R$ whose derivatives of order $m$ are continuous and bounded.

In this paper, we survey some of our work on the near Whitney extension problem \cite{D} in $\mathbb R^n$. Thereafter, we survey some of our work on weighted 
$L_p(\mathbb R),\, 1<p\leq \infty$ convergence of orthonormal expansions in $\mathbb R$ \cite{D1} and present a result of \cite{J}. The motivation for doing this is motivated by 
interesting connections between Whitney extension theorems, Taylor series and Fourier expansions.  Finally, we raise various open questions to study.
\end{abstract}
\medskip


\section{Whitney, Taylor}

\subsection{Whitney extension problem}

Throughout this paper $d\geq 2$, $n,m\geq 1$ and $|.|$ denotes the Euclidean norm on $\mathbb R^n$. Typically, most constants (most often denoted by $c, c_1,..$) are always small enough
and are independent of $d,n$ and possibly other quantities. The context will be clear.
\footnote{By small enough, we mean bounded by a controlled 
small constant. The analogy of large enough is clear.}. If this is not explicitly stated, the context will be clear.
The same letter may denote different constants and/or functions at any given time. The context will be clear. 
The classical Whitney extension problem asks the following: Let $\phi:E\to \mathbb R$ be a map defined on an arbitrary set $E\subset \mathbb R^n$. How can one decide whether $\phi$ extends to a map $\Phi:\mathbb R^n\to \mathbb R$ which agrees with $\phi$ on $E$ and is in $C^m(\mathbb R^n)$, the space of functions from $\mathbb R^n$ to $\mathbb R$ whose derivatives of order $m$ are continuous and bounded.

In order to state one answer to this problem from the paper \cite{W}, we assume that $E$ is closed. We then ask for a theorem that says that it is possible to extend a given real valued function with domain $E$ in such a way as to have prescribed derivatives at the points of $E$. The set $E$ will in general, lack a differentiable structure and to this end we need a careful statement of Taylor's theorem on $\mathbb R^n$ which provides an essential ingredient for the aforementioned theorem's statement and proof. The theorem we state is that of Hassler Whitney \cite{W} and is a partial converse to Taylor's theorem.

\subsection{Taylor to Whitney}

A $n$-dimensional multi-index is an $n$-tuple of non-negative integers $\alpha=(\alpha_1,....,\alpha_n)$. We write $|\alpha|=\alpha_1+...+\alpha_n$,
$\alpha !=\alpha_1 !...\alpha_n !.$ For $x\in \mathbb R^n$, we write $x^{\alpha}=x_1^{\alpha_1}...x_n^{\alpha_n}$.
Given a function $f:\mathbb R^n\to \mathbb R$ in $C^m(\mathbb R^n)$, we set
\[D^{\alpha}f:=\frac{\partial^{|\alpha|}f}{\partial x^{\alpha_1}...\partial x_n^{\alpha_n}},
\] when well defined.
\medskip

We have Taylor's theorem:\, Let $f:\mathbb R^n\to \mathbb R$ in $C^m(\mathbb R^n)$ at a point $a\in \mathbb R^n$.
Then there exist functions $h_{\alpha}:\mathbb R^n\to \mathbb R$ with multi index $\alpha$ such that for $x\in \mathbb R^n$,
\[
f(x)=\sum_{|\alpha|\leq m}\frac{D^{\alpha}f(a)}{\alpha !}(x-a)^{\alpha}+\sum_{|\alpha|=m+1}h_{\alpha}(x)(x-a)^{\alpha}
\] and $\lim_{x\to a}h_{\alpha}(x)=0.$

Differentiating and replacing $h_{\alpha}$ as needed (still denoted by $h_{\alpha}$), we obtain with $f_{\alpha}:=D^{\alpha}f$, for each multi index $\alpha$,
\begin{itemize}
\item[(1)] $f_{\alpha}(x)=\sum_{|\beta|\leq m-|\alpha|}\frac{f_{\alpha+\beta}(a)}{\beta !}(x-a)^{\alpha}+h_{\alpha}(x,a).$
with $h_{\alpha}=o(|x-a|)$ uniformly as $x\to a$.
\end{itemize}

We note that (1) may be regarded as purely a compatibility condition between the functions $f_{\alpha}$ which must be satisfied in order for these functions to be the coefficients of the Taylor series of the function $f$. This important observation allows us to deal with the problem of the closed set $E$, in general lacking 
a differentiable structure in the following Whitney theorem of Hassler Whitney \cite{W}
(below) which follows in part from it.
\medskip

Let $\phi_{\alpha}:E\to \mathbb R$ be a collection of functions defined on a closed set $E\subset \mathbb R^d$ with all multi-indices 
$|\alpha|\leq m$ satisfying the compatibility condition (1) for all $x,a\in \mathbb R^d$. Then there exists a function $\Phi:\mathbb R^d\to \mathbb R$ in 
$C^m(\mathbb R^d)$ such that the following hold:
\begin{itemize}
\item $\Phi=\phi_{0}$.
\item $D^{\alpha}\Phi=\phi_{\alpha}$.
\item $\Phi$ is real-analytic in $\mathbb R^{d}\setminus E$.
\end{itemize}
\medskip

{\bf Fourier series versus Taylor series}
\medskip

A holomorphic function in an anulus containing the unit circle has a Laurent series about zero which generalizes the Taylor series of a holomorphic function in a neighborhood of zero. When restricted to the unit circle, this Laurent series gives a Fourier series of the corresponding periodic function. This explains the connection between the Cauchy integral formula and the integral defining the coefficients of a Fourier series.  The connection between the two provided by the Cauchy integral formula is therefore quite remarkable; one takes an integral of $f$
over the unit circle and it tells you information about the behavior of $f$ at the origin. But this is more a magic property of holomorphic functions than anything else. A variant of this argument would be the following:  Assume that the Taylor expansion of a real valued function $f$:
$f(x)=\sum_{k=0}^{\infty}a_{k}x^k$ is convergent for some $|x|>1$. Then $f$ can be extended in a natural way into the complex plane by writing 
$f(z)=\sum_{k=0}^{\infty}a_kz^k$ with $|z|\leq 1$. Restricting to the unit circle $|z|=1$, define $F(\phi):=f(e^{i\phi})$,  that is, consider $f$ as a function 
of the polar angle $\phi$. $F$ has period $2\pi$ and its Fourier expansion is nothing else but 
$F(\phi)=\sum_{k=0}^{\infty}a_ke^{ik\phi}$ where the $a_k$ are the Taylor coefficients of the real valued map $x\longmapsto f(x)$ 
\bigskip

Our ultimate goal from this paper is to begin an exploration of relations between Witney extensions and orthonormal expansions in line with the above discussion and the interesting paper \cite{J}.
To this end, moving forward we will survey our work in \cite{D1, D} and then state an interesting result of \cite{J} (See Section 5)  and raise some pertinent open questions related to these ideas. 
\medskip

We begin with:

\section{The near distorted Whitney extension problem}
\subsection{The problem}

In this section, we provide a small taiste of some of the work in \cite{D} and refer the reader to the monograph \cite{D} for the full adventure! 

\begin{itemize}
\item A map $A:\mathbb R^{d}\to \mathbb R^{d}$ is an improper Euclidean motion (rigid motion, distance preserving transformation) if there exist $M\in O(d)$ and a translation $x_0\in \mathbb R^d$ so that for every $x\in \mathbb R^d$, $A(x)=Mx+x_0$. 
\item If $M\in SO(d)$, then $A$ is a proper (orientation preserving) Euclidean motion. Here, $O(d)$ and $SO(d)$ are respectively the orthogonal and special orthogonal groups. A Euclidean motion can either be proper or improper. 
\item A map $f:\mathbb R^d\to \mathbb R^d$, will be called a $c$-distortion if there exists $c>0$ small enough depending on $d$ so that 
$(1-c)|x-y|\leq |f(x)-f(y)|\leq |x-y|(1+c)$ for every $x,y\in \mathbb R^d$. Note that $f$ is non-rigid.
\end{itemize}

The problems we are concerned with are: 
\begin{itemize}
\item[(1)] Let $\phi:E\to \mathbb R^d$ be a smooth near distortion \footnote{When we use the term "near distortion", we mean $c$-distortion for small enough $c$.} on a finite set $E\subset \mathbb R^d$ with certain geometry.
How to decide whether $\phi$ extends to a smooth, one to one and onto near distortion $\Phi:\mathbb R^d\to \mathbb R^d$ which agrees with $\phi$ on $E$ and with Euclidean motions in $\mathbb R^d$. 
\item[(2)] Let $U\subset R^n$ be open and let $E\subset U$ be a compact set with certain geometry. Let $\phi:U\to \mathbb R^n$ be a smooth near isometry. How to decide if there exists a smooth one-to-one and onto near isometry
$\Phi:\mathbb R^n\to \mathbb R^n$ which agrees with $\phi$ on $E$ and agrees with Euclidean motions on $\mathbb R^n$.
\end{itemize}

In \cite{D}, we study (1-2) in detail. This is a fascinating problem with connections to 
harmonic analysis, best approximation in algebraic varieties, algebraic geometry in particular, tensor analysis, alignment, clustering, the Procrustes alignment data problem, the manifold data learning problem and BMO (the space of maps of bounded mean variation). See also \cite{A,D1,K9,K10,K12,K13,K14}. 
\medskip

\subsection{Some variants}
Problem (1) below is an example of (1) above: 

\begin{problem}
Let us be given a positive constant $c$ small enough depending on $d$. Does there exist a positive constant $c'$ small enough depending on $c$ so that the following holds.
Given two sets of $k\geq 1$ distinct points in $\mathbb R^d$, $\left\{y_{1},...,y_{k}\right\}$ and $\left\{z_{1},...,z_{k}\right\}$. 
Suppose for every $1\leq i,j\leq k$,
\beq
(1-c')<\frac{|z_{i}-z_{j}|}{|y_{i}-y_{j}|}\leq (1+c').
\label{e:near}
\eeq

\begin{itemize}
\item[(1)] Does there exist a smooth one-to-one and onto $c$-distortion $\Phi:\mathbb R^d\to \mathbb R^d$ 
which obeys $\Phi(y_i)=z_i,\, 1\leq i\leq k$? \footnote{Note that this is an interpolation problem as well a non-rigid alignment data problem.}
\item[(2)] $\Phi$ agrees with Euclidean motions on $\mathbb R^d$.
\item[(3)] Can one say something about how $c, c', k, d$ are related?
\end{itemize}
\label{p:near1}
\end{problem}

\begin{problem}
Let us be given a positive constant $c$ small enough depending on $d$. Does there exist a positve constant $c'$ small enough depending on $c$ so that the following holds.
Given two sets of $k\geq 1$ distinct points in $\mathbb R^d$, $\left\{y_{1},...,y_{k}\right\}$ and $\left\{z_{1},...,z_{k}\right\}$.
Suppose for every $1\leq i,j\leq k$, (\ref{e:near}) holds with $c'$.
\begin{itemize}
\item[(1)] Is it possible to find a Euclidean motion $A$ so that uniformly for every $1\leq i\leq k$, 
\[
(1-c)|y_i-z_i|\leq |A(y_i)-A(z_i)|\leq (1+c)|y_i-z_i|
\]
\item[(2)] Can one say something about how $c, c', k, d$ are related?
\end{itemize}
\label{p:near2}
\end{problem}

In the case of isometric extension, Theorem~\ref{p:near2} is less subtle. It is well known, see for example \cite{WW1} that following holds..
\medskip

Let $\left\{y_{1},...,y_{k}\right\}$ and $\left\{z_{1},...,z_{k}\right\}$ be two collections of $k\geq 1$ distinct points in $\mathbb R^{d}$.
Suppose
\[
|z_{i}-z_{j}|=|y_{i}-y_{j}|,\, 1\leq i,j\leq k.
\]
Then, there exists a Euclidean motion, $A:\mathbb R^{d}\to \mathbb R^{d}$ with
\[
A(y_{i})=z_{i},\, 1\leq i\leq k.
\] 
\footnote{In $\mathbb R^d$, every rigid map is a Euclidean motion. The converse is clear.}

\subsection{Problem 1}

Problem~\ref{p:near1} and Problem~\ref{p:near2} are a new way to view interpolation and non-rigid alignment of data. See \cite{D} for the full details.

\begin{figure}[tb]
\centering\includegraphics[width=0.24\textwidth]{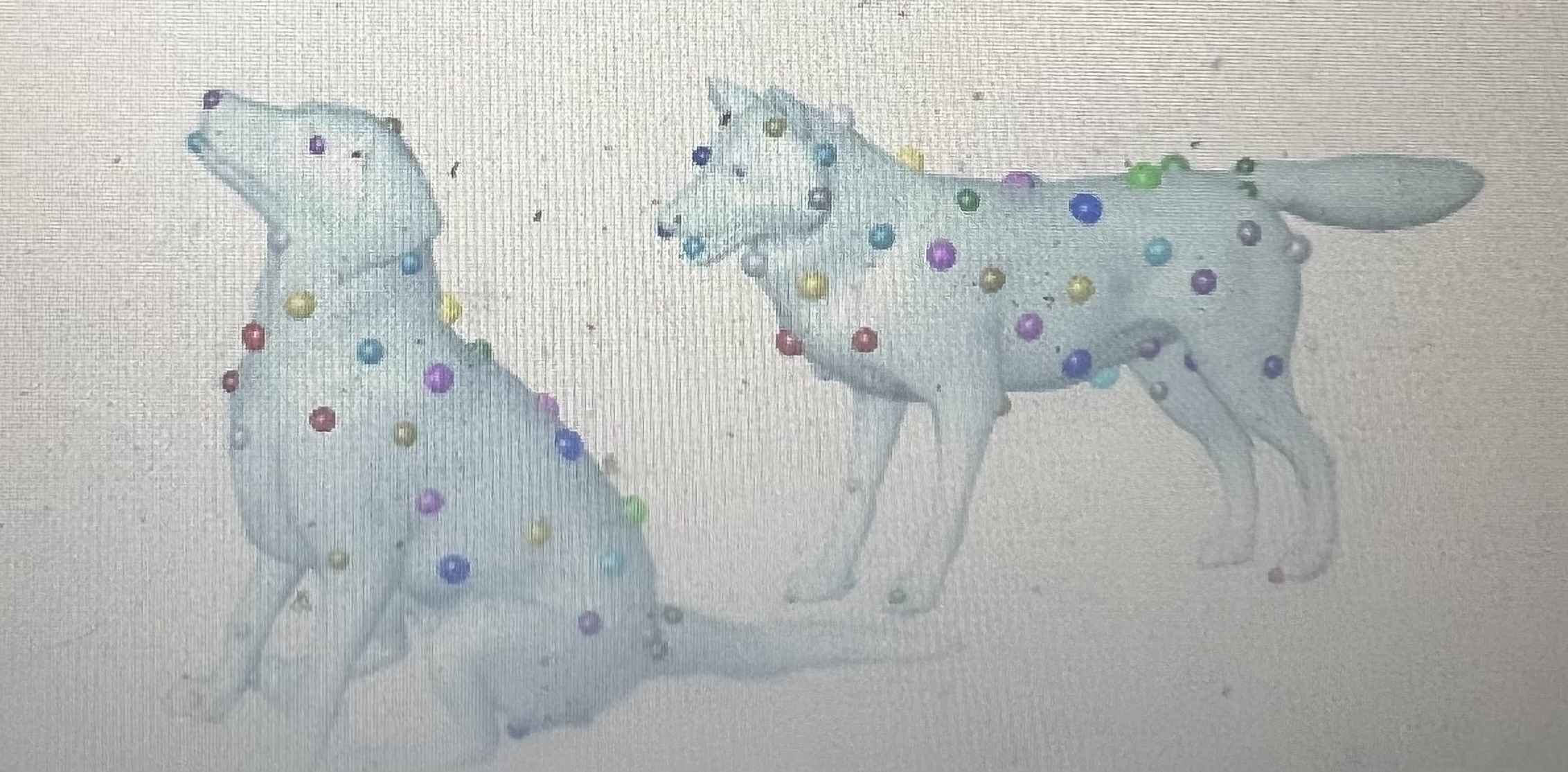}
\caption{Near isometric correspondence shape matching found between a dog and
a wolf. See the reference \cite{LF}}
\end{figure}
\begin{figure}[tb]
\centering\includegraphics[width=0.24\textwidth]{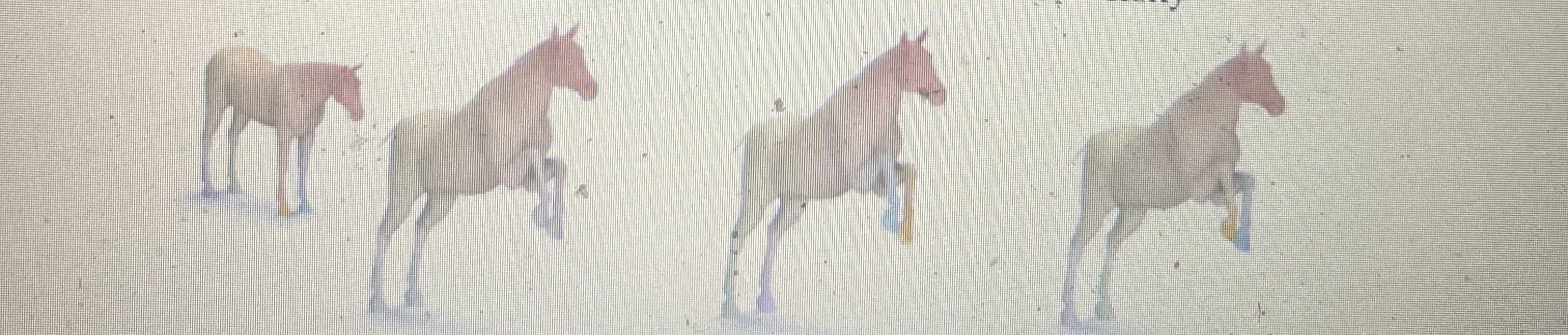}
\caption{Near isometric orrespondences found between horses: See the reference \cite{MO}}
\end{figure}

We now begin to look at Problem~\ref{p:near1} and Problem~\ref{p:near2}. We will translate "smooth" into the idea of a $c$-distorted diffeomorphism.

Given $c$ small enough, a diffeomorphism $f:\mathbb R^{d}\to \mathbb R^{d}$ is a $c$-distorted diffeomorphism if $f$ is a $c$ distortion. \footnote{A diffeomorphism $f:\mathbb R^d\to \mathbb R^d$ is a differentiable one-to-one and onto function with a differentiable inverse.} 

We are now going to provide two examples of $c$-distorted diffeomorphisms.

\subsection{Slow twists}

\begin{exm}
Let $\varepsilon>0$ small enough and $x\in \mathbb R^d$. Let $T_x$ be the block-diagonal matrix given by
\[
\begin{bmatrix}
D_{1}(x) & 0 & 0 & 0 & 0 & 0 \\
0 & D_{2}(x) & 0 & 0 & 0 & 0 \\
0 & 0 & . & 0 & 0 & 0 \\
0 & 0 & 0 & . & 0 & 0 \\
0 & 0 & 0 & 0 & . & 0 \\
0 & 0 & 0 & 0 & 0 & D_{r}(x).
\end{bmatrix}
\]
where for each $i$ either $D_i(x)i$ s the $1\times 1$ identity matrix or else
\[
D_i(x)=\begin{bmatrix}
\cos f_i(|x|) & \sin f_i(|x|) \\
-\sin f_i(|x|) & \cos f_i(|x|)
\end{bmatrix}
\]
for a function $f_i$ of one variable satisfying: \[
A:\,t|f'_{i}(t)|<c\varepsilon,\, t\geq 0.
\]
Here $c$ is small enough*. Then,
\[
f_{ST}(x):=M^{T}T_x(Mx)
\]
for $x\in\mathbb R^d$ and for a fixed $M\in SO(d)$.
is a $\varepsilon$ distorted diffeomorphism. We call $f_{ST}$ a {\it Slow twist} (in analogy to rotations).
\label{e:Example1}
\end{exm}


\subsection{Slides}

\begin{exm}
Let $\varepsilon>0$ be small enough and let $f:\mathbb R^{d}\to \mathbb R$ be a smooth function satisfying the following condition** 
\[
B: |f'(t)|<c\varepsilon,\, t\geq 0.
\]
Here, $c$ is small enough. Consider the function $f_{SL}(t)=t+f(t)$. Then $f_{SL}$ is $\varepsilon$-distorted and we call it a {\it Slide} (in analogy to translations). 
\label{e:Example2}
\end{exm}

Here are two examples: Take a function $f:\mathbb R^2\to \mathbb R$ and suppose (A) holds.
Then 
\[
T_x=\begin{bmatrix}
\cos f(|x|) & \sin f(|x|) \\
-\sin f(|x|) & \cos f (|x|)
\end{bmatrix}.
\]
Take a function $f:\mathbb R^3\to \mathbb R$ and suppose (A) holds. Then 
\[
T_x=\begin{bmatrix}
1 & 0 & 0 \\
0 & \cos(f(|x|)) & \sin(f(|x|))\\
0 & -\sin(f(|x|)) & \cos(f(|x|))\\
\end{bmatrix}.
\]

\subsection{Slow twists: action} 

Here we illustrate the concept of a Slow twist on $\mathbb R^2$. Given a $\varepsilon>0$ and a function $f:\mathbb R^2\to\mathbb R$ so that (A) holds with the function $f$. Define the Slow twist matrix $T_x$ for any $x\in\mathbb R^2$ via
\[ T_x:=\begin{bmatrix}\cos f(|x|) & \sin f(|x|) \\ -\sin f(|x|) & \cos f(|x|)\\ \end{bmatrix}.\]
Then given any pure rotation $M\in\text{SO}(2)$, the following function $f_{St}(x):=M^{T}T_xMx:\mathbb R^2\to\mathbb R^2$ is a Slow twist.

In $2$ dimensions, the rotation $M$ does not affect anything since rotations are commutative on $\R^2$. However, for higher dimensions this is not the case, and hence we leave them in the formulas, but for now we always fix $M$ to be the identity matrix.
For a first set of illustrations, we will look only at one application of a Slow twist with $f$ being an exponential function with differing scaling parameter. 

\begin{figure}[h!]
\centering \includegraphics[scale=0.25]{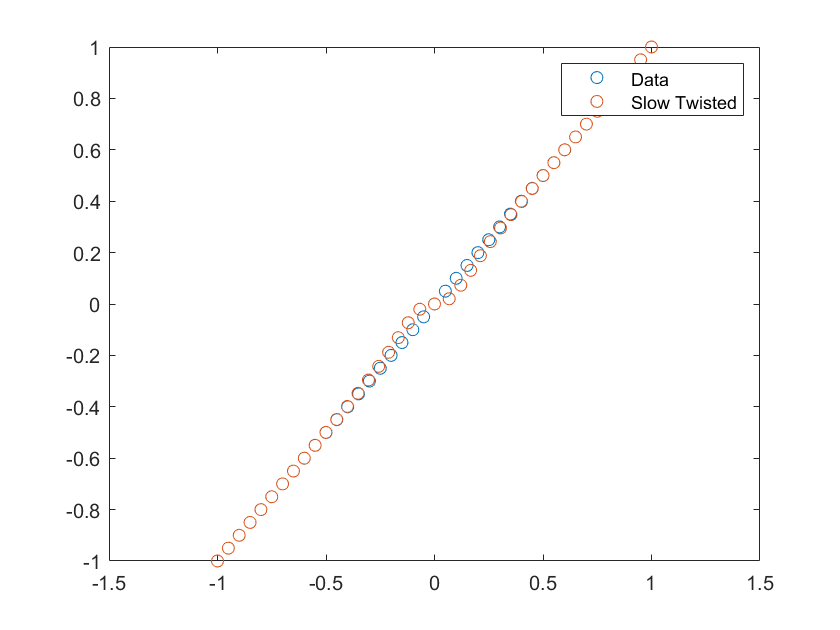} 
\centering \includegraphics[scale=0.25]{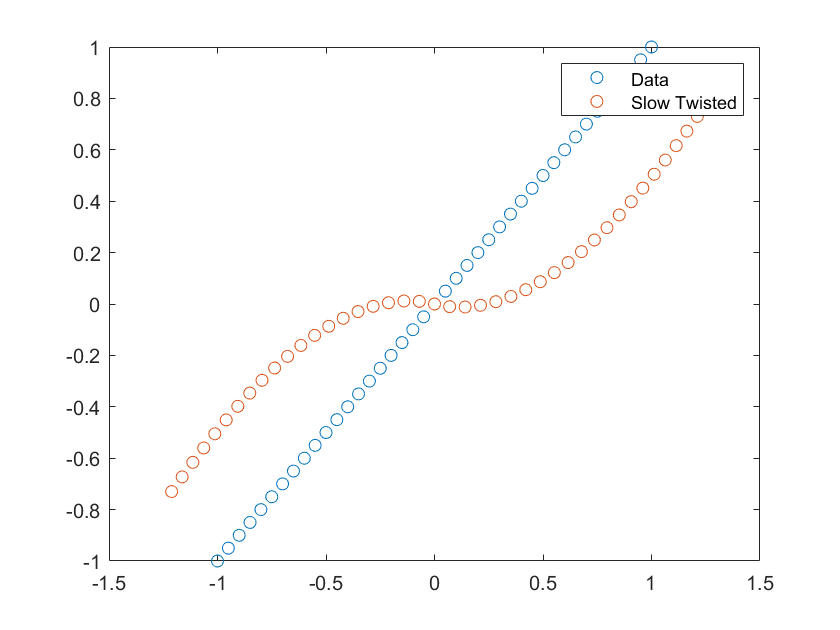} 
\centering \includegraphics[scale=0.25]{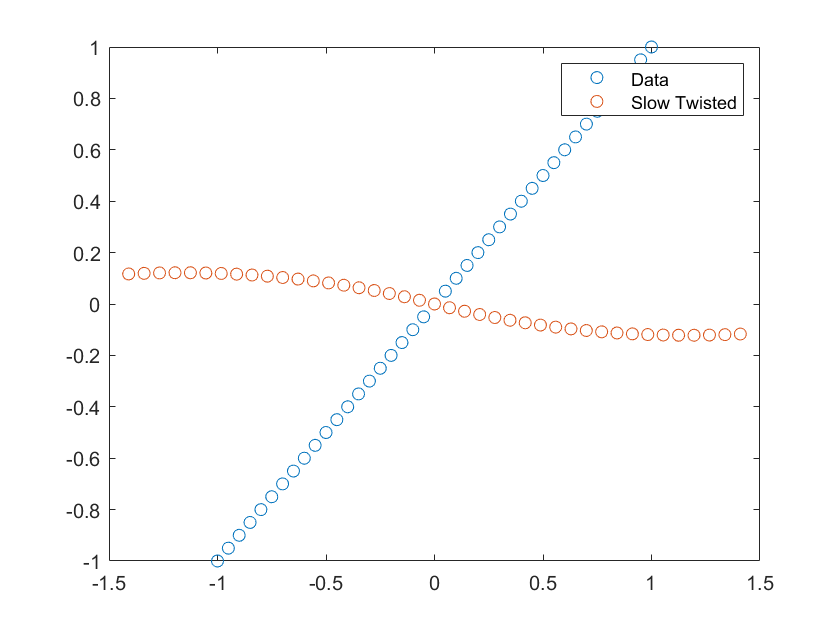}
\caption{Initial points lying on the line $y=x$, and the application of a Slow twist with $f(x)=\exp(-c|x|)$, with $c=10$ (top left), $c=1$ (top middle), and $c=0.1$ top right).}
\end{figure}

For large values of $c$ depending on $d$ the twist is near isometric, and even outside a small enough cube centered at the origin, the points are left essentially fixed. On the other hand, as $c$ tends to $0$, the twist becomes closer to a pure rotation near the origin. Nevertheless, at a far enough distance, the Slow twist $f_{ST}$ will leave the points essentially unchanged. Indeed, the next figures illustrate this:

\begin{figure}[h!]
\centering \includegraphics[scale=0.25]{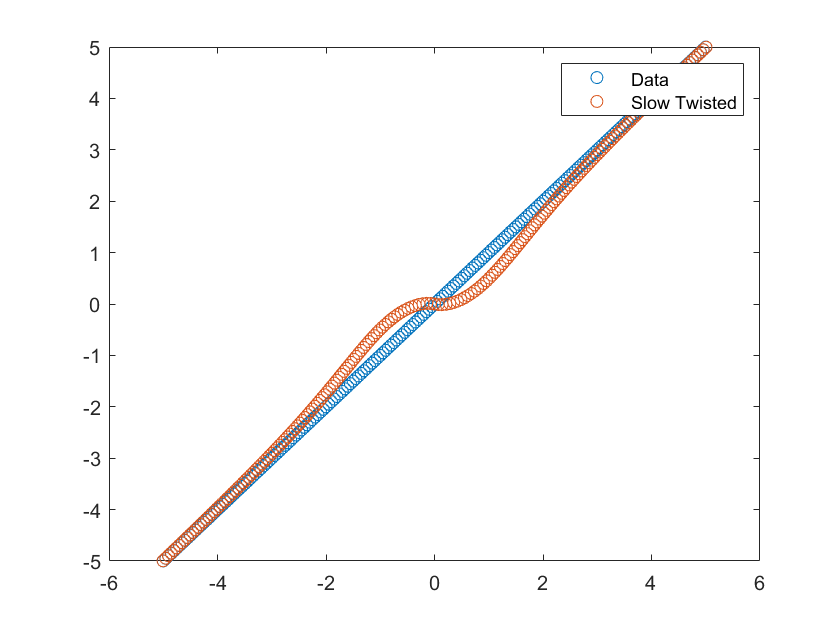} \centering \includegraphics[scale=0.25]{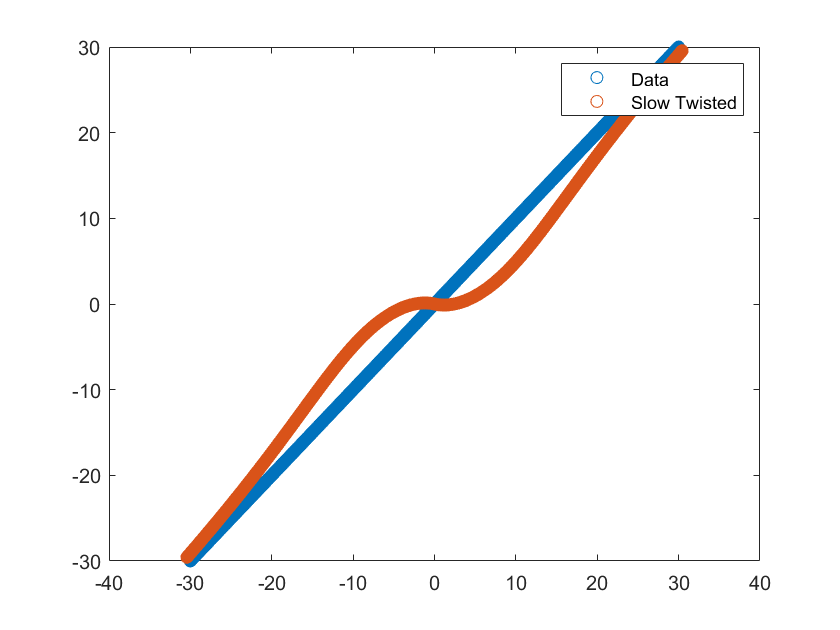}
\caption{Large scale for Slow twists with $f(x)=\exp(-c|x|)$. Left: $c=1$, the twist leaves the points essentially static outside $[-5,5]^2$; Right: $c=0.1$, the twist only starts to leave the points static outside about $[-30,30]^2$.}
\end{figure}

\subsection{Fast twists}

Let us pause to consider what happens when the decay condition on the twist function $f$ (A) is not satisfied; in this case we will dub the twist function a Fast twist for reasons that will become apparent presently.

\begin{figure}[h!]
\centering \includegraphics[scale=0.4]{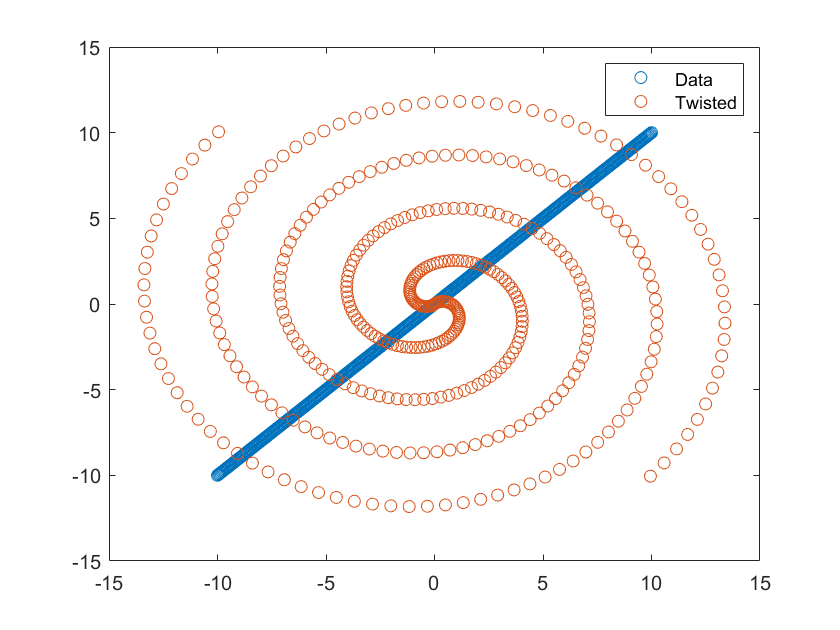}
\caption{Fast twist with function $f(x)=|x|$.}\label{FIG:FastTwist}
\end{figure}

\begin{figure}[h!]
\centering \includegraphics[scale=0.25]{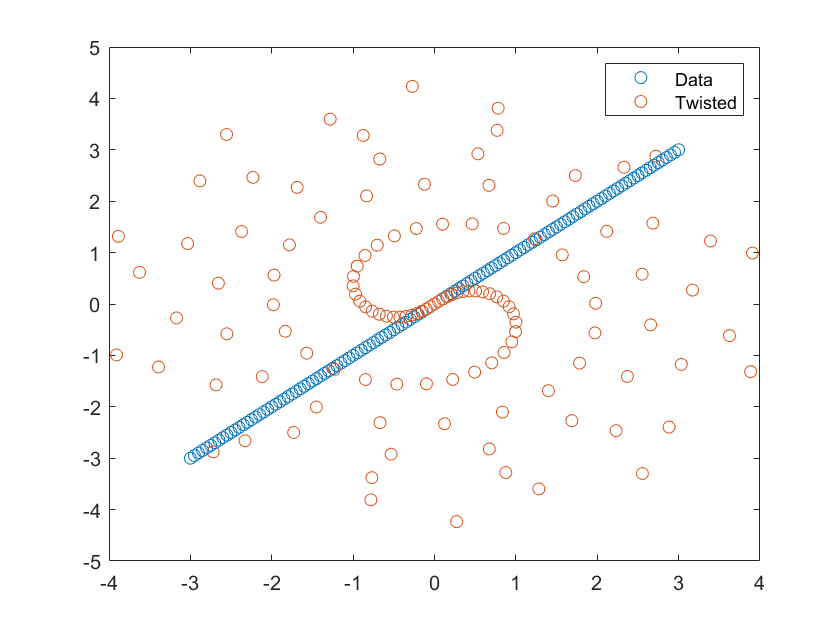} \centering \includegraphics[scale=0.25]{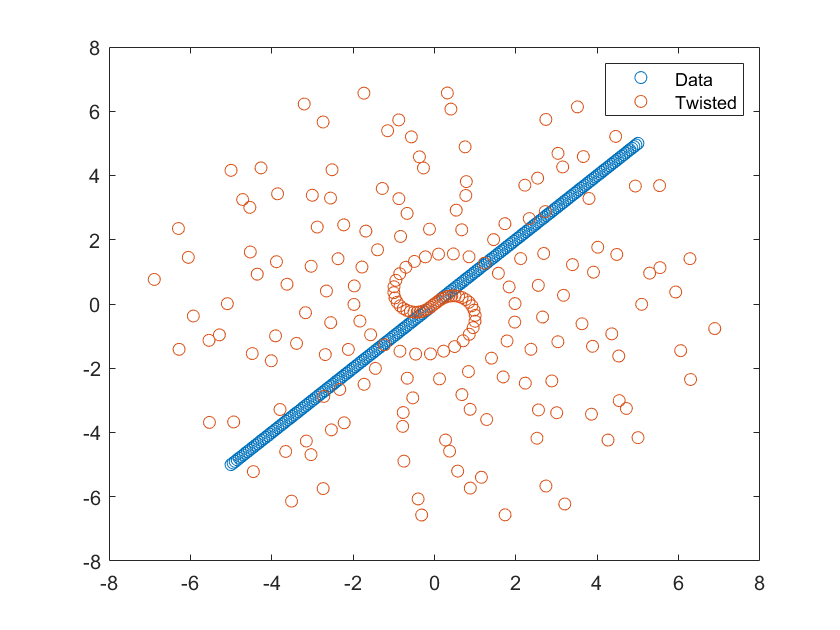} \caption{Fast twist with $f(x)=|x|^2$ for a small enough interval $[-3,3]$ (left) and large interval $[-10,10]$.}\label{FIG:FastTwistSquared}
\end{figure}

One can see that when $f$ is the identity function, the rate of twisting is proportional to the distance away from the origin, and hence there is no way that the twist function will leave points fixed outside of any ball centered at the origin. Likewise, one sees from that the Fast twist with function $f(x)=|x|^2$ rapidly degenerates points into a jumbled mess.

\subsection{Iterated Slow twists}

Here we illustrate what happens when one iteratively applies a Slow twist to a fixed initial point.

\begin{figure}[h!]
\centering \includegraphics[scale=0.5]{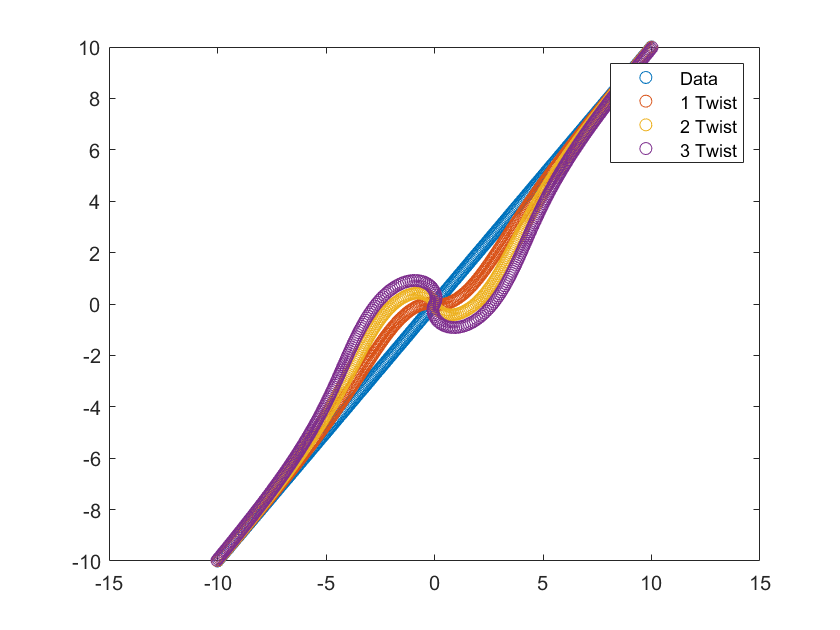} \caption{Iterated Slow twist with $f(x)=\exp(-0.5|x|)$. Shown is the initial points along the line $y=x$, $f_{St}(x)$, $f_{St}\circ f_{St}(x)$, and $f_{St}\circ f_{St}\circ f_{St}(x)$.} \label{FIG:SlowTwistIterated}
\end{figure}

We see an illustration of the fact that the composition of Slow twists remains a Slow twist, but the distortion changes slightly; indeed, notice that as we take more iterations of the exponential Slow twist, we must go farther away from the origin before the new twist leaves the points unchanged.

\begin{figure}
\centering
\animategraphics[loop,autoplay,scale=0.5]{1}{./Twist_}{0}{15}.

\end{figure}

\subsection{Slides: action}

We illustrate some simple examples of Slides on $\R^2$.

\begin{figure}[h!]
\centering
\includegraphics[scale=0.35]{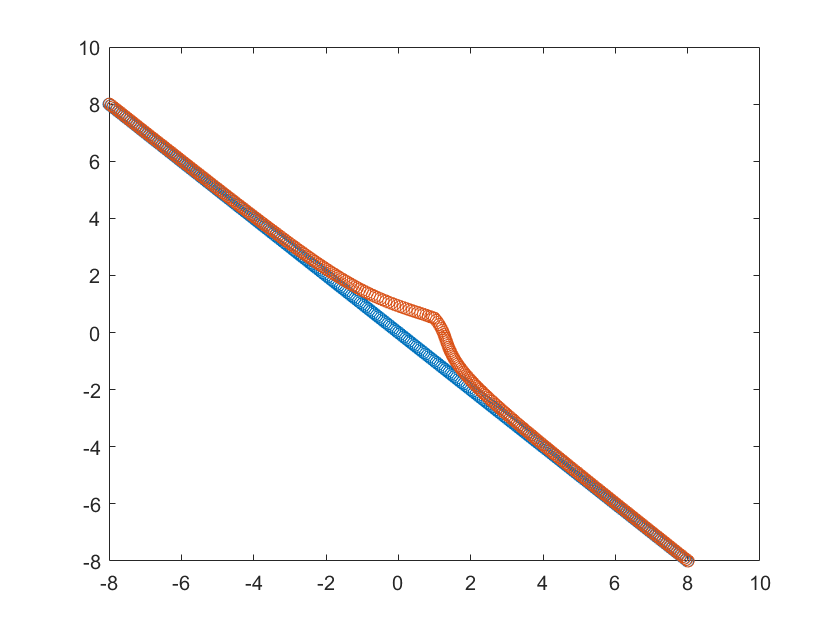}
\caption{Slide with the function $f$ given above.}
\label{FIG:SlideBasic}
\end{figure}

First consider equally spaced points on the line $y=-x$, and the Slide given by the function \[ f(t):=\begin{bmatrix}\frac{1}{1+|t_1|^2}\\ \\ \frac12 e^{-|t_2|} \end{bmatrix}.\]
This is illustrated in Figure \ref{FIG:SlideBasic}.

To give some more sophisticated examples, we consider first the Slide function
\[ f(t) := \begin{bmatrix} e^{-|t_1|} \\ \\ e^{-0.1|t_2|}\end{bmatrix},\] acting iteratively on uniform points along both the lines $y=x$ and $y=-x$.
\begin{figure}[h]
\centering
\includegraphics[scale=0.35]{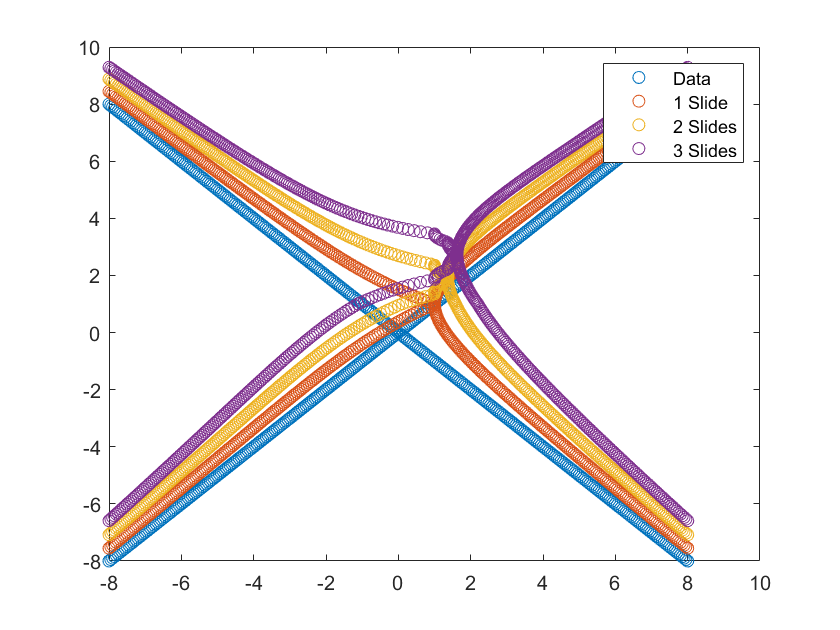}
\caption{Lines $y=x$ and $y=-x$ along with $f_{SL}(x)$, $f_{SL}\circ f_{SL}$ and $f_{SL}\circ f_{SL}\circ f_{SL}(x)$ for $f$.}
\label{FIG:Slide1}
\end{figure}

Similarly, the following figure shows the Slide function
\[ f_2(t) := \begin{bmatrix} 1-e^{-|t_1|} \\ \\ 1-e^{-0.1|t_2|}\end{bmatrix}\]
acting iteratively on uniform points along the lines $y=x$ and $y=-x$.
\begin{figure}
\centering
\includegraphics[scale=0.35]{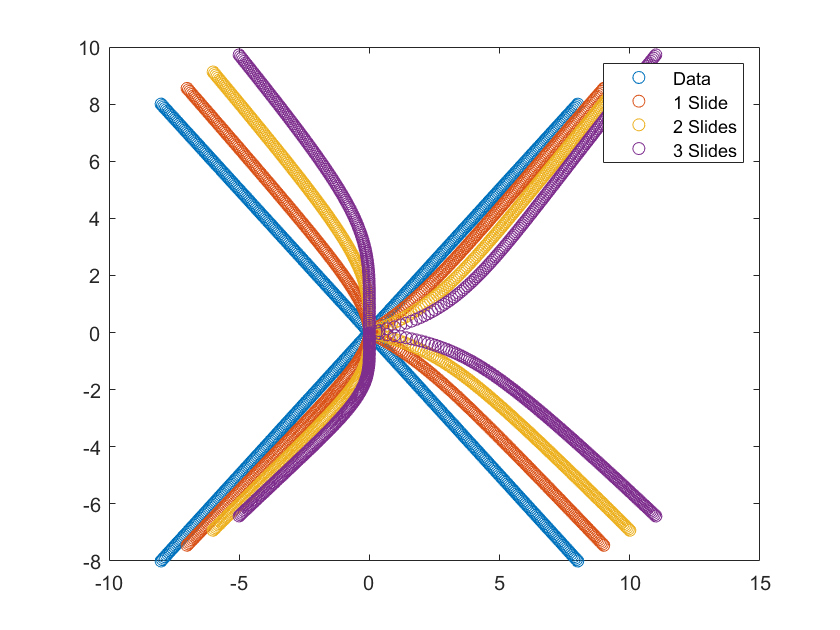}
\caption{Lines $y=x$ and $y=-x$ along with $f_{SL}(x)$, $f_{SL}\circ f_{SL}$ and $f_{SL}\circ f_{SL}\circ f_{SL}(x)$ for $f_2$.}
\label{FIG:Slide2}
\end{figure}

\subsection{Slides at different distances}

To illustrate the effect of the distance of points from the origin, we illustrate here how Slides affect uniform points on circles of different radii.

We use again the asymmetric sliding function 
\[ f(t)=\begin{bmatrix}\frac{1}{1+|t_1|^2}\\ \\ \frac12 e^{-|t_2|} \end{bmatrix}.\]

\begin{figure}
\centering\includegraphics[scale=0.25]{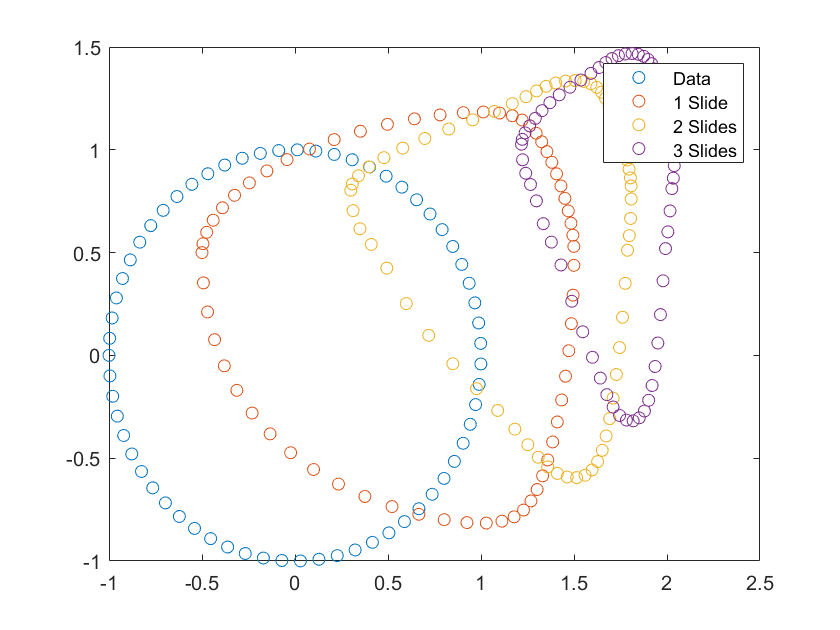} 
\centering\includegraphics[scale=0.25]{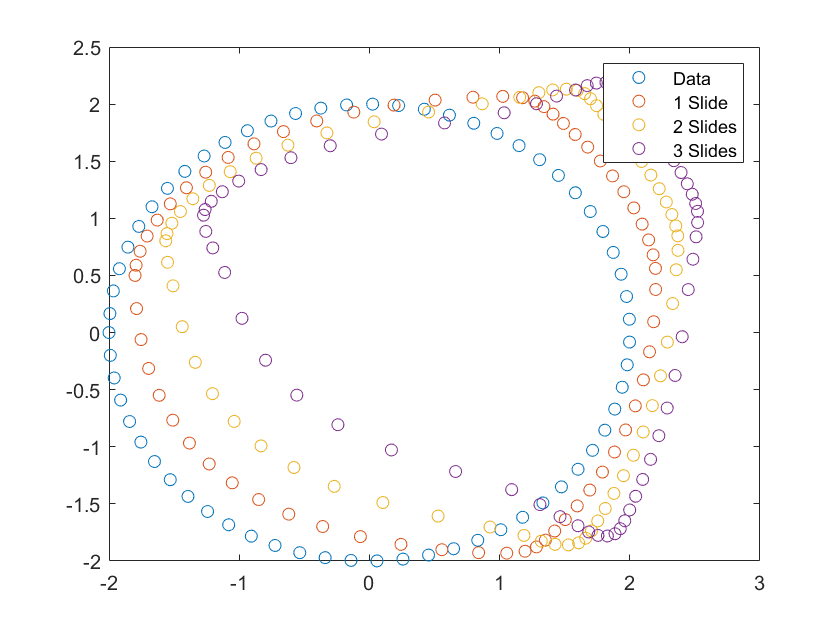} 
\centering \includegraphics[scale=0.25]{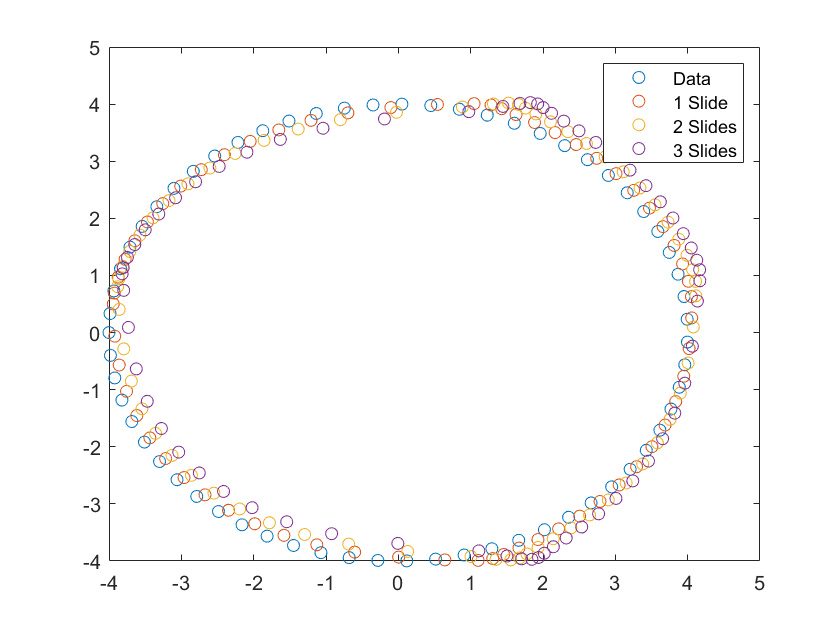}
\caption{Circles under 3 iterated Slides with the function $g$ above, beginning with a circle of radius 1 (top left), 2 (top right) and 4 (bottom).}
\label{FIG:CircleSlides}
\end{figure}

We see that the farther out the points are; i.e., the larger the radius of the initial circle, the less the effect of the Slide, which makes sense given the definition and the fact that the Slides must be $\eps$--distortions of $\R^2$.




\newpage

\subsection{3D motions}

Here we illustrate some of the motions above in $\R^3$.

Here, we construct a generic rotation matrix in $SO(3)$ by specifying parameters $a,b,a_1,d$ satisfying $a^2+b^2+a_{1}^2+d^2=1$, and the rotation matrix $M$ is defined by
\[M = \begin{bmatrix}
a^2+b^2-a_{1}^2-d^2 & 2(ba_{1}-ad) & 2(bd+aa_{1})\\
2(ba_{1}+ad) & a^2-b^2+a_{1}^2-d^2 & 2(a_{1}d- ab)\\
2(bd-aa_{1}) & 2(a_{1}d+ab) & a^2-b^2-a_{1}^2+d^2\\
\end{bmatrix}.\]

As a reminder, our Slow twist on $\R^3$ is thus $M^T_x(Mx)$.

\begin{example}\label{EX:3Dslowtwist.}
Our first example is generated by the rotation matrix $M$ as above with parameters $a=b=\frac{1}{\sqrt{3}}$ and $a_{1}=d=\frac{1}{\sqrt{6}}$, and the Slow twist matrix $St$ as
\[\begin{bmatrix}
1 & 0 & 0 \\
0 & \cos(f(|x|)) & \sin(f(|x|))\\
0 & -\sin(f(|x|)) & \cos(f(|x|))\\
\end{bmatrix},\]
where $f(t) = e^{-\frac{t}{2}}$. Figures \ref{FIG:3DSlowTwist} and \ref{FIG:3DSlowTwistView2} show two views of the twisted motions generated by these parameters.

\begin{figure}[h]
\centering
\includegraphics[scale=0.5]{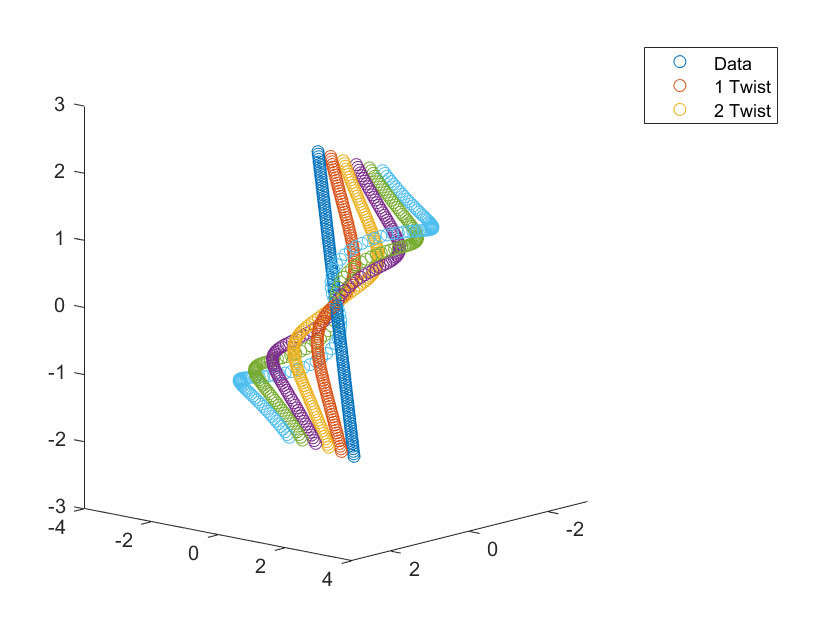}
\caption{Slow twist in $\R^3$.}
\label{FIG:3DSlowTwist}
\end{figure}

\begin{figure}[h]
\centering
\includegraphics[scale=0.5]{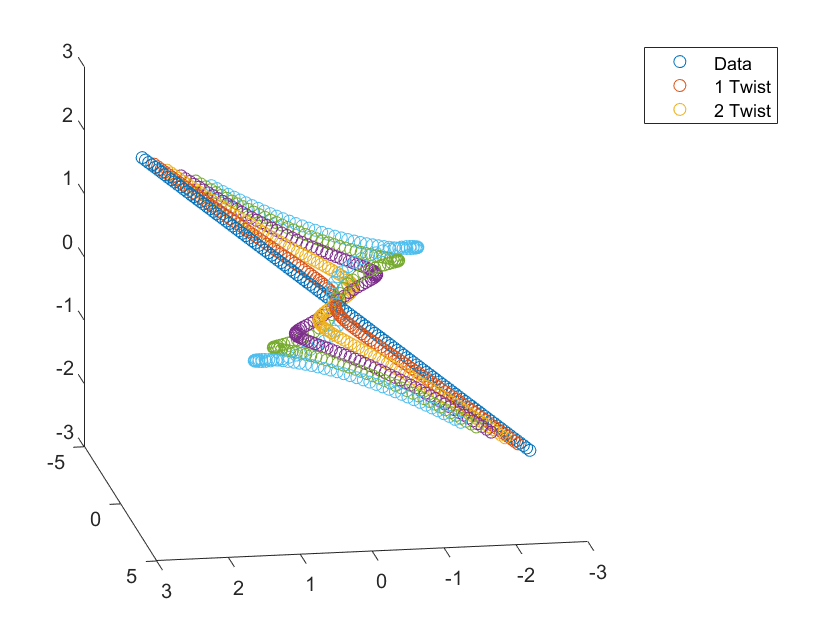}
\caption{Alternate View of Slow twist in $\R^3$.}
\label{FIG:3DSlowTwistView2}
\end{figure}

\end{example}

\begin{figure}
\centering
\animategraphics[loop,autoplay,scale=0.5]{1}{./3DTwist_}{0}{5}
\end{figure}

\subsection{3D Slides}

Here we generate 1000 random points on the unit sphere in $\R^3$ and allow them to move under a Slide formed by 
\[ f(x) = x+\begin{bmatrix}e^{-0.5|x_{1}|}\\ e^{-|x_{2}|}\\ e^{-\frac{3}{2}|x_3|}\\
\end{bmatrix}.\]

\newpage
\newpage

\begin{figure}
\centering
\includegraphics[scale=0.5]{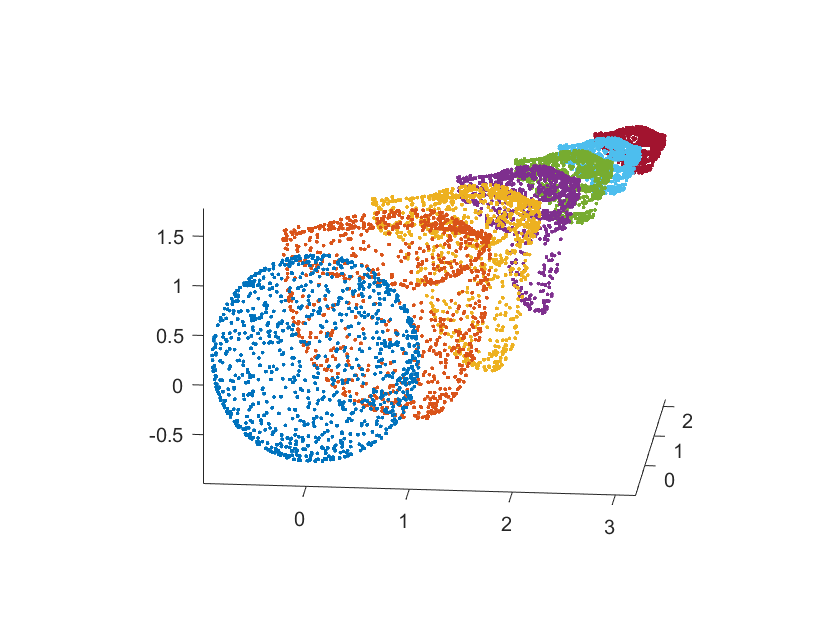}
\caption{Anisotropic Slide on the $2$--sphere.}
\label{FIG:3DSlide}
\end{figure}

From the definition of Slow twists and Slides, the following holds.

\begin{thm}
\begin{itemize}
\item[(1)] Let $\varepsilon>0$ be small enough. There exists $\delta_{1}$ small enough depending on $\varepsilon$ so that the following holds. 
Let $M\in SO(d)$ and let $c_{1}<\delta_{1} c_{2}$. Then there exists a $\varepsilon$- distorted diffeomorphism $f$ with $f(x)=M(x),\, |x|\leq c_{1}$ and 
$f(x)=x,\, |x|\geq c_{2}.$
\item[(2)] Let $\varepsilon>0$. There exists $\delta_{1}$ small enough depending on $\varepsilon$ so that the following holds. 
Let $A(x):=M(x)+x_{0}$ be a proper Euclidean motion and let $c_{3}<\delta_{1} c_{4}$, $|x_{0}|\leq c_5\varepsilon c_{3}$. Then there exists a $\varepsilon$- distorted diffeomorphism $f$ with 
$f(x)=A(x),\, |x|\leq c_{3}$ and $f(x)=x,\, |x|\geq c_{4}$. 
\item[(3)] Let $\varepsilon>0$. There exists $\delta_{1}$ small enough depending on $\varepsilon$ such that the following holds. Let $c_{6}\leq \delta_{1} c_{7}$ and let $x,x'\in \mathbb R^{d}$ with 
$|x-x'|\leq c_{8}\varepsilon c_{6}$ and 
$|x|\leq c_{6}$. Then, there exists a $\varepsilon$-distorted diffeomorphism $f$ with:
\item[(a)] $f(x)=x'$. 
\item[(b)]$f(y)=y,\, |y|\geq c_{7}.$
\end{itemize}
\label{t:mainfd6}
\end{thm}

Moreover, we have:

\begin{thm}
\begin{itemize}
\item Let $E\subset \mathbb R^d$ be a finite set wth ${\rm diameter}(E)\leq 1$. There exists $\varepsilon$-distorted diffeomorphism $\Phi:\mathbb R^d\to \mathbb R^d$ and $c,c'$ depending on $\varepsilon$ with $c$ small enough and $c'$ large enough so that:
\item $\Phi$ coincides with an improper Euclidean motion on $\left\{x\in \mathbb R^d:\, {\rm dist}(x, E)\geq c\right\}$.
\item $\Phi$ coincides with an improper Euclidean motion $A_z$ on $B(z,c')$ for each $z\in E$.
\item $\Phi(z)=z$ for each $z\in E$.
\end{itemize}
\end{thm}

We have then:

\begin{thm}
Let $\varepsilon>0$ be smalll enough.
Then there exist $\delta$ and $\hat{\delta}$ depending on $\varepsilon$ small enough such that the following holds. Let $E\subset \mathbb R^{d}$ be a collection of distinct 
$k\geq 1$ points $E:=\left\{y_{1},...,y_{k}\right\}$. Suppose we are given a function 
$\phi:E\to \mathbb R^{d}$ with
\beq
|x-y||(1+\delta)^{-1}\leq |\phi(x)-\phi(y)|\leq (1+\delta)|x-y|,\, x,y\in E.
\label{e:emotionssa}
\eeq
\begin{itemize}
\item[(1)] If $k\leq d$, there exists a $\varepsilon$-distorted diffeomorphism $\Phi:\mathbb R^{d}\to \mathbb R^{d}$ so that:
\item[(a)] $\Phi$ agrees with $\phi$ on $E$.
\item[(b)] Suppose $y_{i_0}=\phi(y_{i_0})$ for one $i=i_0,\, 1\leq i\leq k$. Then $\Phi(x)=x$, $|x-y_{i_{0}}|\geq \hat{\delta}^{-1/2}{\rm diam}\left\{y_{1},...,y_{k}\right\}$.
\item[(2)] There exists $\delta_{1}$ such that the following holds. Let $E\subset \mathbb R^{d}$ be a collection of distinct 
$k\geq 1$ points $E:=\left\{y_{1},...,y_{k}\right\}$. Suppose that (\ref{e:emotionssa}) holds with $\delta_{1}$. There exists a Euclidean motion $A$ with
\beq
|\phi(x)-A(x)|\leq\varepsilon{\rm diam}(E),\, x\in E. 
\label{e:emotionssaa}
\eeq
If $k\leq d$, then $A$ can be taken as proper.
\end{itemize}
\label{t:mainfd1}
\end{thm}

An important observation in Theorem~\ref{t:mainfd1} is that if $y_{i}=\phi(y_i)$ for one $i=i_0,\, 1\leq i_0\leq k$, the extension $\Phi$ agrees with a Euclidean motion away from the set $E$. What this says is that if the function $\phi$ has a fixed point at one of the points of $E$, then the function $\Phi$ must be essentially rigid away from the set $E$. 
Note the unfortunate restriction of $k\leq d$ in Theorem~\ref{t:mainfd1}(1). The good news is that the problem for this restriction occurs for only some degenerate cases and once removed, an explicit relation between the distortion constants $\varepsilon$ and $\delta$ can be given.
\medskip

In summary: See \cite{D} for more details:

\begin{itemize}
\item Removing the restriction of $k$ in Theorem~\ref{t:mainfd1}(1). Removal of Degenerate Cases. We assume the following:
\item The diameter of the set $E$ is not too large. More precisely: ${\rm diameter}(E)\leq 1$.
\item The points of the set $E$ cannot be too close to each other. More precisely: There exists $0<\eta<1$ so that $|x-y|>\eta$ for every distinct $x,y\in E$.
\item The points of the set $E$ are close to a hyperplane in $\mathbb R^d$. More precisely: 
For a set of $l+1$ points in $\mathbb R^{d}$, with $l\leq d$, say $z_{0},…,z_{l}$ we define
$V_{l}(z_{0},..., z_{l}):={\rm vol}_{l\leq d}({\rm simplex}_{l})$ where ${\rm simplex}_{l}$ is the $l$-simplex with vertices at the points $z_{0},...,z_{l}$.
Thus $V_{l}(z_{0},..., z_{l})$ is the $l\leq d$-dimensional volume of the $l$-simplex with vertices at the points $z_{0},...,z_{l}$.
We write $V_{l}(E)$ as the maximum of $V_{l}(z_0,...,z_l)$ over all 
points $z_{0},z_{1},...,z_{l}$ in $E$. Assume the following: There exists $0<\eta<1$ so that $V_{l}(E)<\eta$.
\item Another interpretation: What is required is that on any $d+ 1$ of the $k$ points which form vertices of a relatively
voluminous simplex, the mapping $\phi$ is orientation preserving. 
\item We then have: $\delta=\exp\left(-\frac{c}{\varepsilon}\right)$.
\end{itemize}

Regarding (2) in Theorem~\ref{t:mainfd1} , it is shown in \cite{D} using tools in approximation theory and algebraic geometry that if the points $E$ are forced lie on ellipse of radius 1, then one can take $\delta=c\varepsilon^c$.
\medskip

We now work with compact sets in open sets in $\mathbb R^n$ and achieve the optimal $\delta=c\varepsilon$.
\medskip

We work with a $C^1$ map $\phi: U\to \mathbb R^n$ with $U\subset \mathbb R^n$ be an open set. For a compact set $E$ in $\mathbb R^n$, 
and $x \in \bbR^{n}$, we write for convenience $d(x):={\rm dist}(x, E)$. Let $\ep>0$.

We assume the following:
\begin{itemize}
\item [(a)] Geometry of $E$: For certain positive constants $c_{0}$, $C_{1}$, $c_{2}$ depending on $n$, the following holds:
Let $x \in \bbR^{n} \backslash E$. If $d(x) \leq c_{0} \, {\rm diam}(E)$, then there exists a ball $B(z,r) \subset E$ such that $| z-x|\leq c_{1} \, d(x)$ and $r \geq c_{2} \, d(x)$.
\item[(b)] Geometry of $\phi$: For $x,y\in E$, $|x-y|(1-\ep)\leq |\phi(x)-\phi(y)|\leq (1+\ep)|x-y|.$
\item[(c)] Controlled constants: $\ep$ is smaller than a small positive constant determined by $c_i$, $i=0,1,2$ and $n$.
\end{itemize}


Here is the main result from \cite{D}.

\begin{thm}
\label{t:0.1}
Under the above assumptions, there exists a $C^{1}$ function $\Phi: \bbR^{n} \to \bbR^{n}$ and a Euclidean motion $A: \bbR^{n} \to \bbR^{n}$, with the following properties:

\begin{itemize}
\item \label{i:0.4} \textit{$(1-c'\ep)\ml x-y \mr \leq \ml \Phi(x) - \Phi(y) \mr \leq (1+c'\ep)\ml x-y \mr$ for all $x,y \in \bbR^{n}$}. $c'$ determined by $c_i,\, i=1,2,3$ and $n$.

\item \label{i:0.5} \textit{$\Phi = \phi$ in a neighborhood of $E$}.

\item \label{i:0.6} \textit{$\Phi = A$ outside $\bl x \in \bbR^{n}: d(x) < c_{0} \br$}.

\item \label{i:0.7} \textit{$\Phi: \bbR^{n} \to \bbR^{n}$ is one-to-one and onto}.

\item \label{i:0.8} \textit{If $\phi \in C^{m}(U)$ for some given $m\geq 1$, then $\Phi \in C^{m}\(\bbR^{n}\)$}.

\item \label{i:0.9} \textit{If $\phi \in C^{\infty}(U)$, then $\Phi \in C^{\infty}(\bbR^{n})$}.
\end{itemize}
\label{t:section6exten}
\end{thm}

Moving forward we are now going to survey some of our work on weighted 
$L_p(\mathbb R),\, 1<p\leq \infty$ convergence of orthonormal expansions in $\mathbb R$, see \cite{D1}. Our survey 
is motivated by our discussion above under the title "Fourier series versus Taylor series", the work of \cite{J} (orthonormal expansions)  and the work on Taylor series (Section 1.2) both of which imply Whitney extension theorems.

\section{Orthonormal expansions on the line $\mathbb R$}
For a positive measure $d\alpha$ on the real line with infinitely many points its support and with all power moments $\int x^l\alpha(t),\, l=0, 1, 2,...$ finite, the Gram-Schmidt process then generates orthonormal polynomials
$p_n(x)=p_n(d\alpha)(x)=\gamma_nx^n+...+\gamma_n,$ satisfying $\int p_n(t)p_l(t)d\alpha(t)=\delta_{ln}$. \footnote{Recall from the beginning that we always have $n\geq 1$.} Here,
\[
\delta_{l,n}:=\left\{
\begin{array}{ll}
0, & l\neq n \\
1, & l=n
\end{array} \right.
\]
Here we set $p_{-1}=0$, $p_0=\left(\int d\alpha(t)dt \right)^{-1/2}$.
The orthonormal polynomials satisfy the recurrence.
\[
xp_{n}(x)=\alpha_np_{n+1}(x)+\alpha_{n-1}p_{n-1}(x),\, x\in \mathbb R
\]
where $\alpha_n:=\frac{\gamma_{n-1}}{\gamma_n}$, $\gamma_{0}=1$ and ${\rm lim}_{n\to\infty}\frac{\alpha_{n+1}}{\alpha_n}=1.$ Now for functions $f:\mathbb R\to \mathbb R$ for which $f(t)t^l,\, l\geq 0$ is integrable with respect to $d\alpha$, we can form the formal 
orthornomal polynomial expansion $S(f):=S(f,d\alpha)=\sum_{l=0}^{\infty}c_lp_l$ where for $1\geq 0,\, c_l:=\int fp_ld\alpha$ and by 
$S_n(f):=S_n(f,d\alpha)=\sum_{j=0}^{n-1}c_jp_j$, the $(n-1)$th partial sum of $S(f)$. 
We are interested in studying how to understand for fixed $f$, convergence of $S_n(f)$ to $S(f)$ as $n\to \infty$.

\section{Translation to weighted $L_p(\mathbb R),\, 1<p\leq \infty$ convergence of $S_n(f)$ to $S(f)$ as $n\to \infty$ for fixed admissible $f$.}

We set $d\alpha(t)=w(t)^2dt$ with $w$ a positive even function on $\mathbb R$ of smooth polynomial decay for large argument. 
For $1<p\leq\infty$ and fixed weight $w$, we will always work with functions $f$ so that the norm below is finite. We will call such functions admissible.
\[
|fw|_{p(\mathbb R)}:=
\left\{
\begin{array}{ll}
{\rm sup}_{x\in \mathbb R}|fw(x)|, & p=\infty \\
\left(\int_{\mathbb R}|(fw)(x)|^pdx\right)^{1/p}, & 1<p<\infty.
\end{array}
\right.
\]
Given $u,x\in \mathbb R$, we set
\[
\Delta_{u}(f)(x):=f(x+u)-f(x)
\]
and for a fixed $\gamma\in
\mathbb R$, we will also define
\[
u_{\gamma}(y):=(1+|y|)^{\gamma},\, y\in \mathbb R.
\] 
\subsection{A class of admissible weights.}

We work with the following class of admissible weights. See \cite{L4, L5, D1}.
\begin{dfn}
A weight function $w={\rm exp}(-Q):\mathbb R\to (0,\infty)$
will be called {\it admissible}
if each of the following conditions below is satisfied:
\begin{itemize}
\item[(a)] $Q:={\rm log}(1/w)$
is continuously differentiable, even and satisfies $Q(0)=0$;
\item[(b)] $Q'$ is nondecreasing in $\mathbb R$ with
\[
\lim_{x\to \infty}Q(x)=\lim_{x\to -\infty}Q(x)=\infty.
\]
Assume that there exists $\eta>1$ with
\[
\eta<\frac{xQ'(x)}{Q(x)}\leq C,\, x\in \mathbb R\backslash\left\{0\right\}.
\]
\item[(c)] For every $\varepsilon>0$, there exists $\delta>0$ such that
for every $x\in \mathbb R\backslash\{0\}$,
\[
\int_{x-\delta |x|}^{x+\delta |x|}
\frac{|Q'(s)-Q'(x)|}{|s-x|^{3/2}}ds
\leq \varepsilon|Q'(x)|.
\]
\end{itemize}
\label{d:admiss}
\end{dfn}

Definition~(\ref{d:admiss}), defines a general class of
even weights for which our
results hold. \footnote{Later, we will look at a paper \cite{J} which shows a deep connection between Whitney extensions and orthogonal expansions for Laguerre orthogonal polynomials.}
The
weak regularity and
smoothness conditions on $w$ above are needed, for bounds on $p_n$, and its
zeroes and are used heavily in our proofs of the results below taken from \cite{D1}.
Note that Definition~(\ref{d:admiss}) does not require $Q''$ to exist.
Instead, we require only a mild local Lipshitz $1/2$
condition on $Q'$.
Note that $Q$ grows as a polynomial for large argument and $w$ is of smooth polynomial decay for large argument. 
We note, as an easily absorbed example, that
$w_{\beta}(x):={\rm exp}\left(-|x|^{\beta}\right),\,\beta>1,\, x\in \mathbb R$
is an example of an admissible weight. The case $\beta=2$ is just the Hermite weight. See \cite{L4}. 
\medskip

\subsection{The numbers $a_u$ and $\alpha_n$}
In studying weighted polynomial approximation for admissible weights $w$ on the
line, an important
role is played by the scaled endpoints $\pm a_u$
of the support of the equilibrium measure for $w^2$ (see for example \cite{D1,L4,L5})
and the asymptotic behavior of the quotient $\frac{\alpha_n}{a_n}$. More precisely:
\medskip

Given $u>0$, we define the real number $a_u$ by the positive root of the
equation
\[
u=\frac{2}{\pi}\int_0^1\frac{a_utQ'(a_ut)}{\sqrt{1-t^2}}dt,\, u>0.
\]
It is known, see \cite{L4, L5, D1}, that $a_u$ is
uniquely defined, strictly increasing in $(0,\infty)$ with
\[
\lim_{u\to\infty}a_u=\infty 
\]
and of polynomial growth for large argument. For example, for the
weight, $w_{\beta}(x):={\rm exp}\left(-|x|^{\beta}\right)$ above , it is known that
\[
a_u\sim u^{1/\beta}.
\]
One of the important properties of the number $a_u$ which allows the analysis the convergence of $S_n(f)$ to $S(f)$ in $L_{p}(\mathbb R)$ is the fact that, see
\cite{L4, L5, D1}, that
\[
|Pw|_{\infty(\mathbb R)}=|Pw|_{\infty[-a_n,a_n]}
\]
and
\[
|(Pw)(x)|_{\infty(|x|\geq sa_n)}\leq \exp(-cn)|Pw|_{\infty[-a_n,a_n]}
\]
for every fixed $s>1$ and for every polynomial $P$ of degree
at most $n$ with $1<p<\infty$ analogues. By differentiation, it not difficult to see why the above should hold for ${\rm exp}\left(-|x|^{\beta}\right)$. Indeed, these 
"Infinite-finite" inequalities allow a shift in 
various weighted polynomial approximations on $\mathbb R$ to approximations on sequences of compact intervals. It is also known, see \cite{L4, L5, D1} that
$\lim_{n\to\infty}\frac{\alpha_n}{a_n}=1/2$.

\subsection{Main Results in \cite{D1}}
Following are the main results of \cite{D1}. 

\subsubsection{Necessary results: $L_{\infty}(\mathbb R)$}

\begin{thm}
Let $w$ be an admissible weight and
$B,b\in \mathbb R$ with $b<B$. Then for
\begin{eqnarray}
&& {\rm sup}_{n\in\Omega}|S_n(f)wu_b(x)| \\
\nonumber && \leq c|fwu_B|_{\infty(\mathbb R)}
\label{e:ness1}
\end{eqnarray}
to hold for some infinite subsequence $\Omega\subseteq \mathbb N$, for all
$x\in \mathbb R$ and admissible $f$ for which the right-hand side of (4.4) is finite
it is necessary that
\begin{equation}
B>0
\label{e:ness2}
\end{equation}
and uniformly
\begin{equation}
a_n^{b-{\rm min}\left\{B,1\right\}}n^{1/6}C_{B,n}=O(1)
\label{e:ness3}
\end{equation}
where
\[
C_{B,n}:=\left\{
\begin{array}{ll}
1, &B\neq 1 \\
\log n, &B=1.
\end{array}\right.
\]
\label{t:thm1}
\end{thm}

\subsubsection{Sufficiency results: $L_{\infty}(\mathbb R)$}

\begin{thm}
Let $w$ be an admissible weight, $b\leq 0$ and
assume (4.5) and (4.6). Let $x\in \mathbb R$ and assume moreover that
\begin{equation}
\frac{\alpha_{n+1}}{\alpha_n} = 1 +
O\left( \frac{1}{n}\right),\qquad n\to\infty.
\end{equation}
\footnote{(4.7) and (4.9) are known for classes of smooth admissible weights and the weight $w_{2n}$, see \cite{D20, Deift}}. Then there exists an 
infinite subsequence $\Omega\subseteq \mathbb N$
so that for $n\in \Omega$
\begin{eqnarray}
&& |(S_n(f)wu_b)(x)|\\
\nonumber && \leq c
\left[{\rm log}n|fwu_{B}|_{\infty(\mathbb R)}+
\int_{-1}^{1} \left|\frac{[wu_b\Delta_y(f)](x)}{y}dy\right|\right]
\end{eqnarray}
for admissible $f$ for which the right-hand side of (4.8) is finite. If in addition,
(4.9) holds below,
\begin{equation}
\frac{\alpha_n}{a_n} = \frac12 \left[1 +
O\left( \frac{1}{n^{2/3}} \right) \right],\qquad n\to\infty,
\end{equation}
then for $n\geq 1$
\begin{eqnarray}
&& |(S_n(f)wu_b)(x)|\\
\nonumber && \leq c
\left[{\rm log}n|fwu_{B}|_{\infty(\mathbb R)}+
\int_{-1}^{1}\left|
\frac{[wu_b\Delta_y(f)](x)}{y}dy\right|\right]
\end{eqnarray}
for all admissible $f$
for which (4.10) is finite.
\end{thm}

\begin{thm}
Let $w$ be admissible, $x\in\mathbb R$, $b\leq 0$ and
assume (4.5)-(4.7). Then there
exists an infinite subsequence $\Omega\subseteq\mathbb N$ such that for $n\in \Omega$
\begin{equation}
|S_n[f]wu_b(x)|
\leq c\left[{\rm log}n|fwu_{B}|_{\infty(\mathbb R)}+|f'wu_b|_{\infty(\mathbb R)}
\right]
\end{equation}
for admissible $f:\mathbb R\to\mathbb R$ for which the right side of (4.11) is finite. Moreover, if (4.9) holds in addition, then for $n\geq 1$,
\begin{equation}
|S_n[f]wu_b(x)|
\leq c\left[{\rm log}n|fwu_{B}|_{\infty(\mathbb R)}+|f'wu_b|_{\infty(\mathbb R)}
\right]
\end{equation}
for all admissible functions $f:\mathbb R\to\mathbb R$ for
which the right-hand side of (4.12) holds.
\end{thm}

Note that (4.8) and (4.11) are for a subsequence whereas (4.10) and (4.12) is for all $n\geq 1$. 
\medskip

Now we look at $L_p$ analogues of Theorems (4.6-4.7).

\subsection{$L_{p}(\mathbb R)$ results}

\begin{thm}
Let $w$ be admissible, $b,B\in \mathbb R$ with $b\leq B$,
$1<p<\infty$ and $C_{B,n}$
as defined in (4.6). 
Then for
\begin{equation}
{\rm sup}_{n\in \Omega}|S_{n}[f]wu_b|_{p(\mathbb R)}\leq
c|fwu_B|_{p(\mathbb R)}
\end{equation}
to hold for some infinite subsequence $\Omega\subseteq \mathbb N$ and for all admissible $f$ for which (4.13) is finite 
it is necessary that
\begin{equation}
b<1-1/p,\, B>-1/p.
\end{equation}
In addition:
\begin{itemize}
\item[(a)] If $p<4/3$ then necessarily,
\[
a_n^{{\rm max}\left\{b,-1/p\right\}-B}n^{1/6(4/p-3)}C_{B,n}=O(1).
\]
\item[(b)] If $p=4/3$ or $4$ then necessarily $b<B$.
\item[(c)] If $p>4$ then necessarily
\[
a_n^{b-{\rm min}\left\{B,1-1/p\right\}}n^{1/6(1-4/p)}C_{B,n}=O(1).
\]
\end{itemize}
Moreover, (4.7), (4.14), (a-c) are sufficient to give (4.13) for all admissible $f$ for which the right-hand side of (4.13) is finite. In addition, if (4.9) also holds then 
\begin{equation}
{\rm sup}_{n\geq 1}|S_{n}[f]wu_b|_{p(\mathbb R)}\leq
c|fwu_B|_{p(\mathbb R)}.
\end{equation}
and
\begin{equation}
\lim_{n\to \infty}|(S_{n}[f]-f)wu_b|_{p(\mathbb R)}=0
\end{equation}
for all continuous $f:\mathbb R\to\mathbb R$ satisfying
\begin{equation}
\lim_{|x|\to\infty}|fwu_{B+\delta}(x)|=0
\end{equation}
for some $\delta>1$.
\end{thm}
\section{Orthonormal expansions to Whitney}
In this section, we illustrate an interesting connection between certain Whitney extensions and certain Laguerre polynomial orthonormal expansions detailed in the paper \cite{J}. By $\mathbb R_{+}^d$, we mean the usual cartesian product space $(0,\infty)^d$, by 
$\overline{\mathbb R_{+}^d}$ its closure $[0,\infty)^d$ and similary define ${\mathbb N}_0^d$ 
with ${\mathbb N}_0:={\mathbb N}\cup \left\{0\right\}$. $S(\mathbb R_{+}^d)$ consists of all real valued functions $f\in C^{\infty}(\mathbb R_{+}^d)$ such that all 
derivatives $D^p f,\, p\in {\mathbb N}_0^d$ extend to continuous functions on $\overline{\mathbb R_{+}^d}$ and satisfying
${\rm sup}_{x\in {\mathbb R}_{+}^d}x^k|D^p f(x)|<\infty,\, \forall k,p\in {\mathbb N}_0^d.$ The space $S(\mathbb R^d)$ is defined similarly.
The main result in the paper \cite{J} deals with the existence of a Whitney extension of the space $S(\mathbb R_{+}^d)$ onto the space $S(\mathbb R^d)$ .
A crucial tool in dealing with this problem are certain Laguerre polynomial orthonormal expansions. 
\medskip

We recall for $n=0,1,2,...$, the functions $P_n(x):=\frac{e^x}{n!}\left(\frac{d}{dx}\right)^{n}(e^{-x}x^n),\, x>0$ are the Laguerre polynomials, 
${\cal P}_n(x):=P_n(x)e^{x/2}$ are the Laguerre functions and the set $\left\{{\cal P}_n(x),\, n=0,1,...\right\}$ is an orthonormal basis 
for $L^2(0,\infty)$, the space of squared integrable functions on $(0,\infty)$. \cite{S}.
\medskip

In \cite{J}, the following is proved: Convergence of the Laguerre series in the space $S(\mathbb R_{+}^d)$ and this is then used as a tool to obtain the existence of a Whitney extension of the space $S(\mathbb R_{+}^d)$ onto the space $S(\mathbb R^d)$.

\begin{itemize}
\item[(1)] Convergence of the Laguerre series in the space $S(\mathbb R_{+}^d)$: For $f\in S(\mathbb R_{+}^d)$, let $a_n(f):=\int_{\mathbb R_{+}^d}
f(x){\cal P}_{n}(x)dx.$ Then $f(x)=\sum_{n\in {\mathbb N}_0^d}a_n(f){\cal P}_n(x)$ converges absolutely for all $x\in S(\mathbb R_{+}^d)$. 
\item[(2)] Existence of a Whitney extension of the space $S(\mathbb R_{+}^d)$ onto the space $S(\mathbb R^d)$: The restriction 
mapping $f\to f|_{\mathbb R_{+}^d}:\,S(\mathbb R^d)\to S(\mathbb R_{+}^d)$ is a topological onto homomorphism. Moreover, 
the space $S(\mathbb R_{+}^d)$ is topologically isomorphic to the quotient space $S(\mathbb R^d)/N$ where
\[
N:=\left\{f\in S(\mathbb R^d)\,|\,{\rm supp}\,f\subseteq (\mathbb R^d\setminus \mathbb R_{+}^d)\right\}.
\]
\end{itemize}

\section{Open Questions}

\begin{itemize}
\item[(a)] It would be of interest to develop a framework for the ideas in Section (5) or the near Whitney extension problem using the work in Section (4) replacing Laguerre polynomials by the weighted orthogonal polynomials and subsequent orthonormal expansion theory there.
\item[(b)] It would also be of interest to do the same in (a) except using weighted orthogonal polynomials with respect to different fast decreasing weights over different domains, for example 
the even Erdos type weights over $\mathbb R$, of faster than smooth polynomial decay with large argument, the even Pollazek weights which are of faster than smooth polynomial decay
close to the endpoints of $(-1,1)$, classes of exponential (possibly non-symmetric weights) on real intervals $(a,b)$ with $a<0<b$ and other classical orthogonal polynomials 
on $\mathbb R^n$. See \cite{D2,L4,L,L1}
\end{itemize}

\end{document}